\documentclass[journal,twoside]{IEEEtran}
\usepackage{cite}

\ifCLASSINFOpdf
\else
\fi
\usepackage[cmex10]{amsmath}
\usepackage{array}

\usepackage{mdwmath}
\usepackage{mdwtab}
\usepackage[caption=false,font=footnotesize]{subfig}
\usepackage{fixltx2e}

\usepackage{stfloats}
\usepackage{url}

\usepackage[american,fulldiode]{circuitikz} 
\usepackage{graphicx} 
\usepackage{multirow}
\usepackage{placeins}
\usepackage{amsfonts}
\usepackage{bbm}
\usepackage{algorithm}
\usepackage[noend]{algpseudocode}
\usepackage{multirow}

\usepackage{slashbox}
\usepackage{amssymb}

\usepackage{lipsum,amsmath,multicol}

\newcommand{\matpower}{M{\sc atpower}}
\usepackage{multirow}

\makeatletter
\def\ps@headings{%
\def\@oddhead{\mbox{}\scriptsize\rightmark \hfil \thepage}%
\def\@evenhead{\scriptsize\thepage \hfil \leftmark\mbox{}}%
\def\@oddfoot{}%
\def\@evenfoot{}}
\def\BState{\State\hskip-\ALG@thistlm}
\makeatother

\usepackage{varwidth}
\newcolumntype{M}[1]{>{\begin{varwidth}[t]{#1}}l<{\end{varwidth}}}

\begin{document}
%
\title{Sparsity-Exploiting Moment-Based Relaxations of the Optimal Power Flow Problem}
%
%
%

\author{Daniel K. Molzahn, \IEEEmembership{Member, IEEE} and~Ian A. Hiskens, \IEEEmembership{Fellow, IEEE}
\thanks{University of Michigan Department of Electrical Engineering and Computer Science: molzahn@umich.edu, hiskens@umich.edu}
}
\maketitle

\begin{abstract}
Convex relaxations of non-convex optimal power flow (OPF) problems have recently attracted significant interest. While existing relaxations globally solve many OPF problems, there are practical problems for which existing relaxations fail to yield physically meaningful solutions. This paper applies moment relaxations to solve many of these OPF problems. The moment relaxations are developed from the Lasserre hierarchy for solving generalized moment problems. Increasing the relaxation order in this hierarchy results in ``tighter'' relaxations at the computational cost of larger semidefinite programs. Low-order moment relaxations are capable of globally solving many small OPF problems for which existing relaxations fail. By exploiting sparsity and only applying the higher-order relaxation to specific buses, global solutions to larger problems are computationally tractable through the use of an iterative algorithm informed by a heuristic for choosing where to apply the higher-order constraints. With standard semidefinite programming solvers, the algorithm globally solves many test systems with up to 300 buses for which the existing semidefinite relaxation fails to yield globally optimal solutions.
\end{abstract}



%
\IEEEpeerreviewmaketitle

\section{Introduction}
\label{l:introduction}
%
%
%
%

\IEEEPARstart{T}{he} optimal power flow (OPF) problem determines an optimal operating point for an electric power system in terms of a specified objective function (typically generation cost per unit time), subject to both network equality constraints (i.e., the power flow equations, which model the relationship between voltages and power injections) and engineering limits (e.g., inequality constraints on voltage magnitudes, active and reactive power generations, and line flows). While the OPF problem is often augmented with security constraints that ensure robustness to contingencies (see, e.g.,~\cite{Capitanescu,stott2012,phan2014,platbrood2014}), the formulation considered in this research does not consider contingencies.

The OPF problem is generally non-convex due to the non-linear power flow equations~\cite{bernie_opfconvexity} and may have local solutions~\cite{bukhsh_tps}. Non-convexity of the OPF problem has made solution techniques an ongoing research topic. Many OPF solution techniques have been proposed, including successive quadratic programs, Lagrangian relaxation, genetic algorithms, particle swarm optimization, and interior point methods~\cite{opf_survey,opf_litreview1993I,opf_litreview1993II,matpower,frank2012opf,ferc4}. Some of these techniques are quite mature and capable of finding at-least-locally optimal solutions to many large-scale OPF problems with reasonable computational burden. (For instance,~\cite{phan2014} and~\cite{platbrood2014} report computationally tractable methods for finding at-least-locally-optimal solutions to security-constrained OPF problems with approximately 3,000 and 9,000 buses, respectively.)

However, while typical local solution techniques often in fact find globally optimal solutions~\cite{molzahn_lesieutre_demarco-global_optimality_condition,phan2014}, they may fail to converge or converge to a local optimum. See, for instance, the example problems and discussion in~\cite{bukhsh_tps} as well as the five-bus system in~\cite{bernie_opfconvexity}. For many of these problems, \matpower's~\cite{matpower} default interior point solver with default options and a variety of initialization heuristics either fails to converge or finds locally optimal solutions. See also~\cite{ferc5}, which reports on a study by the Federal Energy Regulatory Commission of convergence characteristics for a variety of commercial solvers, OPF problem formulations, and initialization heuristics.

Recently, significant attention has focused on a semidefinite relaxation of the OPF problem~\cite{lavaei_tps}. If the relaxed problem satisfies a rank condition (i.e., the relaxation is said to be ``exact'' or ``tight''), the global solution to the original OPF problem can be determined in polynomial time. Prior OPF solution methods do not guarantee finding a global solution in polynomial time. Further, infeasibility of a relaxation certifies infeasibility of the OPF problem, which is a capability not available with typical existing solution techniques. Additionally, unlike local solution techniques whose convergence characteristics generally depend on a chosen initialization, the semidefinite relaxation provides the global solution regardless of the choice of initialization when the relaxation is exact. While not as mature as existing solution techniques, semidefinite programming approaches thus have substantial advantages over traditional techniques.

However, the rank condition is not satisfied for all practical OPF problems~\cite{bukhsh_tps,hicss2014}. For such problems, the relaxation provides a lower bound on the optimal objective value but does not provide physically meaningful decision variables (i.e., voltage phasors). The bounds obtained from the semidefinite relaxation are often close to the global optimum and are therefore useful for many applications (e.g., calculating voltage stability margins~\cite{molzahn_lesieutre_demarco-pfcondition} and determining the potential suboptimality of an OPF solution that is only guaranteed to be locally optimal). However, determining both the globally optimal objective value and the globally optimal decision variables is important in many contexts. 

To address problems for which the rank condition is not satisfied, this paper presents \emph{moment relaxations}\footnote{The terminology \emph{moment relaxation}, adopted from~\cite{lasserre2001,lasserre_book}, refers to the relaxation's derivation from a special case of the \emph{generalized moment problem}.} that globally solve a broader class of OPF problems than existing relaxations. Building on the results of~\cite{pscc2014} (many of which are independently studied in~\cite{patrick} and~\cite{ibm_paper}), this paper describes a method for globally solving large OPF problems by exploiting sparsity and only applying computationally intensive ``higher-order'' moment relaxations to specific buses.

Much of the related work in this area focuses on sufficient conditions for which existing convex relaxations are exact~\cite{low_tutorial1,low_tutorial2,madani2014}. While the sufficient conditions developed thus far are promising, they only apply to a limited subset of problems. 

For more general cases,~\cite{bitar_allerton2013} proposes a heuristic method for finding a global optimum that is ``hidden'' in a higher-rank subspace of solutions to the semidefinite relaxation. However, the semidefinite relaxation yields an optimal objective value strictly less than the global minimum of some OPF problems~\cite{hicss2014}. For such cases, other heuristics may obtain at-least-locally optimal solutions~\cite{madani2014,bitar_allerton2013} with the optimal objective value of the semidefinite relaxation indicating the potential suboptimality.

While deserving of further study, heuristics eliminate the global optimality guarantee that is one of the main advantages of the semidefinite relaxation. We propose an alternative moment relaxation that, when exact, yields the global optimum. Using polynomial optimization theory~\cite{lasserre2001,lasserre_book}, moment relaxations globally solve a broad class of OPF problems, including many problems for which existing relaxations are not exact. Moment relaxations exploit the fact that the OPF problem is composed of polynomials in the voltage phasor components and is therefore a polynomial optimization problem. 

Global solution of a broader class of OPF problems has a computational cost. Whereas the matrix in the semidefinite relaxation of~\cite{lavaei_tps} is composed of all second-order combinations of the voltage phasor components, the moment relaxation's matrices are composed of higher-degree combinations. The semidefinite program for the order-$\gamma$ moment relaxation of an $n$-bus system has a positive semidefinite constraint on a $k\times k$ matrix, where $k = \left(2n+\gamma\right)! / \left( \left(2n\right)! \gamma!\right)$ (i.e., this matrix is composed of all combinations of voltage components up to order $2\gamma$). For example, the moment matrices for the \mbox{first-,} second-, and third-order moment relaxations of a 10-bus system have size $21 \times 21$, $231 \times 231$, and $1771 \times 1771$, respectively, as compared to $20\times 20$ for the semidefinite relaxation of~\cite{lavaei_tps}. Thus, the computational requirements of the moment relaxations can be substantially larger than the semidefinite relaxation of~\cite{lavaei_tps}, especially for high-order relaxations.

Fortunately, experience with small systems suggests that low-order relaxations globally solve a broad class of OPF problems, including problems for which the semidefinite relaxation of~\cite{lavaei_tps} is not exact. As an example of the effectiveness of the moment relaxations, consider the 9-bus OPF problem in~\cite{bukhsh_tps}. \matpower~\cite{matpower} with the default interior point solver and default solver options either fails to converge or converges to one of three local optima depending on the initialization.\footnote{Initialization heuristics included 1.)~a ``flat start'' with unity voltage magnitudes and zero voltage angles, 2.)~the solution to the linear ``DC'' OPF approximation, 3.) a power flow solution calculated using active power injections at the midpoints of the generators' operating ranges, 4.)~a power flow solution calculated using power injections corresponding to an economic dispatch, and 5.)~a power flow solution calculated using the power injections resulting from a DC OPF.} (The local optima have objective values that are 10.0\%, 37.5\%, and 38.1\% greater than the global optimum.) The semidefinite relaxation of~\cite{lavaei_tps} yields a lower bound that is 11\% less than the global optimum. Thus, an existing convex relaxation and a typical interior point technique both perform poorly for this problem while a second-order relaxation finds the global solution~\cite{pscc2014}. The capabilities of low-order relaxations for small OPF problems are further described in~\cite{pscc2014}, which includes an exploration of the feasible spaces of second-order relaxations, and independently in~\cite{patrick} and~\cite{ibm_paper}.

However, large OPF problems are computationally intractable even for low-order relaxations. Solving existing semidefinite relaxations of large OPF problems requires exploiting power system sparsity. Using a matrix completion decomposition, existing semidefinite relaxations are computationally tractable for problems with thousands of buses~\cite{jabr11,molzahn_holzer_lesieutre_demarco-large_scale_sdp_opf}. Na\"{\i}ve application of related techniques for the moment relaxations~\cite{waki2006} enables solution of systems with up to approximately forty buses. Solving larger systems requires more judicious use of the higher-order relaxation. Exploiting the observation that power injection ``mismatches'' typically occur only in small regions of realistic large OPF problems~\cite{hicss2014}, this paper applies a higher-order relaxation to specific buses. This enables global solution of large problems. With standard semidefinite programming solvers, the proposed approach is successfully applied to OPF problems with up to 300 buses for which the semidefinite relaxation of~\cite{lavaei_tps} fails to yield the globally optimal decision values. Further improvements in solving larger OPF problems may be achieved by combining emerging semidefinite programming solvers with the method proposed in this paper.

Selective application of the moment relaxation is independently proposed in~\cite{ibm_paper}. The method described in~\cite{ibm_paper} is limited to second-order relaxations of OPF problems with less than 40 buses due to a computationally expensive subproblem and lack of concurrent exploitation of sparsity.

After introducing the OPF problem formulation in Section~\ref{l:opf_formulation}, we describe the moment relaxations in Section~\ref{l:msdp_overview}. Section~\ref{l:msdp_sparse} then presents the method for globally solving large OPF problems by exploiting sparsity and only applying the higher-order relaxations to specific buses. Section~\ref{l:results} presents results from the proposed method. Section~\ref{l:conclusion} concludes the paper and discusses future research directions.

\section{OPF Problem Formulation}
\label{l:opf_formulation}

We first present an OPF formulation in terms of rectangular voltage coordinates, active and reactive power generation, and apparent-power line-flow limits. Consider an $n$-bus power system, where $\mathcal{N} = \left\lbrace 1, 2, \ldots, n \right\rbrace$ is the set of all buses, $\mathcal{G}$ is the set of generator buses, and $\mathcal{L}$ is the set of all lines. Let $P_{Dk} + j Q_{Dk}$ represent the active and reactive load demand at each bus~$k \in \mathcal{N}$. Let $V_k = V_{dk} + j V_{qk}$ denote the voltage phasors in rectangular coordinates at each bus~$k \in \mathcal{N}$. Superscripts ``max'' and ``min'' denote specified upper and lower limits. Buses without generators have maximum and minimum generation set to zero. Let $\mathbf{Y} = \mathbf{G} + j \mathbf{B}$ denote the network admittance matrix. Shunt conductances and susceptances at bus~$k$ contribute to the diagonal element $\mathbf{Y}_{kk}$.


\pagebreak
The power flow equations describe the network physics:

\begin{subequations}
\small
\begin{align}\nonumber
P_{Gk} = & f_{Pk}\left(V_d,V_q\right) = V_{dk} \sum_{i=1}^n \left( \mathbf{G}_{ik} V_{di} - \mathbf{B}_{ik} V_{qi} \right) &  &  \\ 
\label{opf_Pbalance}  & + V_{qk} \sum_{i=1}^n \left( \mathbf{B}_{ik}V_{di} + \mathbf{G}_{ik}V_{qi} \right) + P_{Dk}  \\ \nonumber
Q_{Gk} = & f_{Qk}\left(V_d,V_q \right) = V_{dk} \sum_{i=1}^n \left( -\mathbf{B}_{ik}V_{di} - \mathbf{G}_{ik} V_{qi}\right) \\
\label{opf_Qbalance} & + V_{qk} \sum_{i=1}^n \left( \mathbf{G}_{ik} V_{di} - \mathbf{B}_{ik} V_{qi}\right) + Q_{Dk}
\end{align}
\end{subequations}

Define a convex quadratic cost of active power generation:

\begin{equation}\label{objfunction}
f_{Ck}\left(V_d,V_q\right) = c_{k2} \left(f_{Pk}\left(V_d,V_q\right)\right)^2 + c_{k1} f_{Pk}\left(V_d,V_q\right) + c_{k0}
\end{equation}

\noindent Note that while we focus on minimization of a quadratic function of active power generation, one can substitute other cost functions (e.g., loss minimization, voltage regulation, convex piecewise-linear generation cost functions, reactive power dispatch, etc.) for~\eqref{objfunction}. (The moment relaxation approach described in this paper is applicable for any polynomial or convex piecewise-polynomial objective function.)

Define a function for squared voltage magnitude:

\begin{equation} \label{opf_Vsq}
\left(V_{k}\right)^2 = f_{Vk}\left(V_d, V_q\right) = V_{dk}^2 + V_{qk}^2
\end{equation}

Squared apparent-power line flows are polynomial functions of the voltage components $V_d$ and $V_q$. To account for flow limits on transformers with non-zero phase shifts and/or off-nominal voltage ratios, we model the line from bus~$l$ to bus~$m$ as a $\Pi$-model circuit with series admittance $g_{lm} + j b_{lm}$ and total shunt admittance $g_{sh,lm} + j b_{sh,lm}$ in series with an ideal transformer with a specified complex turns ratio $1\colon \tau_{lm} e^{j\theta_{lm}}$ as in~\cite{matpower}. (Note that the conductance $g_{sh,lm}$ in the $\Pi$-model is generally neglected in typical power system data sets, and that shunt susceptances $b_{sh,lm}$ are often neglected for transformers.)


\begin{subequations}
\small
\begin{align}
\nonumber  & P_{lm} = f_{Plm}\left(V_d,V_q\right) = \left( V_{dl}^2 + V_{ql}^2\right) \left(g_{lm} + \frac{g_{sh,lm}}{2}\right) / \tau_{lm}^2 \\ \nonumber
& \quad +
\left(V_{dl}V_{dm} + V_{ql}V_{qm}\right) \left(b_{lm}\sin\left(\theta_{lm} \right) - g_{lm}\cos\left(\theta_{lm} \right) \right) / \tau_{lm} \\
\label{Plm}& \quad +
\left(V_{dl}V_{qm} - V_{ql}V_{dm}\right) \left(g_{lm}\sin\left(\theta_{lm}\right) + b_{lm}\cos\left(\theta_{lm}\right)\right) / \tau_{lm}
\end{align}
\begin{align}
\nonumber & P_{ml} = f_{Pml}\left(V_d,V_q\right) = \left(V_{dm}^2 + V_{qm}^2 \right)\left(g_{lm}  + \frac{g_{sh,lm}}{2} \right) \\ \nonumber & \quad -
\left(V_{dl}V_{dm} + V_{ql}V_{qm} \right)\left(g_{lm}\cos\left(\theta_{lm}\right) + b_{lm}\sin\left(\theta_{lm}\right) \right) / \tau_{lm} \\ \label{Pml}
& \quad + 
\left(V_{dl}V_{qm} - V_{ql}V_{dm} \right) \left(g_{lm}\sin\left(\theta_{lm}\right) - b_{lm}\cos\left(\theta_{lm}\right) \right) / \tau_{lm}
\end{align}
\begin{align}
\nonumber & Q_{lm} = f_{Qlm}\left(V_d,V_q\right) = -\left( V_{dl}^2 + V_{ql}^2\right) \left(b_{lm} + \frac{b_{sh,lm}}{2}\right) / \tau_{lm}^2 \\ \nonumber 
& \quad +
\left(V_{dl}V_{dm} + V_{ql}V_{qm}\right) \left(b_{lm}\cos\left(\theta_{lm} \right) + g_{lm}\sin\left(\theta_{lm} \right) \right) / \tau_{lm} \\ \label{Qlm} & \quad +
\left(V_{dl}V_{qm} - V_{ql}V_{dm}\right) \left(g_{lm}\cos\left(\theta_{lm}\right) - b_{lm}\sin\left(\theta_{lm}\right)\right) / \tau_{lm}
\end{align}
\begin{align}
\nonumber & Q_{ml} = f_{Qml}\left(V_d,V_q\right) =  -\left( V_{dm}^2 + V_{qm}^2\right) \left(b_{lm} + \frac{b_{sh,lm}}{2}\right) \\ \nonumber & \quad +
\left(V_{dl}V_{dm} + V_{ql}V_{qm}\right) \left(b_{lm}\cos\left(\theta_{lm} \right) - g_{lm}\sin\left(\theta_{lm} \right) \right)  / \tau_{lm} + \\ \label{Qml} & \quad + 
\left(-V_{dl}V_{qm} + V_{ql}V_{dm}\right) \left(g_{lm}\cos\left(\theta_{lm}\right) + b_{lm}\sin\left(\theta_{lm}\right)\right) / \tau_{lm} \\ 
\label{Slm} & \left(S_{lm}\right)^2 = f_{Slm}\left(V_d,V_q\right) = \left(f_{Plm}\left(V_d,V_q\right)\right)^2 + \left(f_{Qlm}\left(V_d,V_q\right)\right)^2 \\
\label{Sml} & \left(S_{ml}\right)^2 = f_{Sml}\left(V_d,V_q\right) = \left(f_{Pml}\left(V_d,V_q\right)\right)^2 + \left(f_{Qml}\left(V_d,V_q\right)\right)^2
\end{align}
\end{subequations}

The classical OPF problem is then

\begin{subequations}
\label{opf}
\small
\begin{align}
\label{opf_obj} & \min_{V_d,V_q} \sum_{k \in \mathcal{G}} f_{Ck}\left(V_d,V_q\right) \qquad \mathrm{subject\; to} \hspace{-20pt} & \\
\label{opf_P} &  P_{Gk}^{\mathrm{min}} \leq f_{Pk}\left(V_d,V_q \right) \leq P_{Gk}^{\mathrm{max}} & \forall k \in \mathcal{N} \\
\label{opf_Q} &  Q_{Gk}^{\mathrm{min}} \leq f_{Qk}\left(V_d,V_q \right) \leq Q_{Gk}^{\mathrm{max}} &  \forall k \in \mathcal{N} \\
\label{opf_V} &  \left(V_{k}^{\mathrm{min}}\right)^2 \leq f_{Vk}\left(V_d, V_q\right) \leq \left(\vphantom{V_{k}^{\mathrm{min}}} V_{k}^{\mathrm{max}}\right)^2 &  \forall k \in \mathcal{N}  \\
\label{opf_Slm} &  f_{Slm}\left(V_d,V_q\right) \leq \left(S_{lm}^{\mathrm{max}}\right)^2 &  \forall \left(l,m\right) \in \mathcal{L} \\
\label{opf_Sml} &  f_{Sml}\left(V_d,V_q\right) \leq \left(S_{lm}^{\mathrm{max}}\right)^2 &  \forall \left(l,m\right) \in \mathcal{L} \\
\label{opf_Vref} & V_{q1} = 0
\end{align}
\end{subequations}

\noindent Constraints~\eqref{opf_Slm} and~\eqref{opf_Sml} limit the apparent-power flow at each line terminal. Constraint~\eqref{opf_Vref} sets the reference bus angle to zero.

Note that the OPF problem is often extended to consider contingency, voltage-stability, and transient stability constraints. (See~\cite{Capitanescu,stott2012,phan2014,platbrood2014} for discussions of these and other extensions.) As a starting point for the methods developed in this paper, we only consider OPF problems without these constraints. Future work includes extension to more general OPF problem formulations. (See~\cite{lavaei2011acc} and~\cite{maria2013powertech} for initial research on applications of convex relaxations to OPF problems with contingency constraints.)

Further, the OPF formulation studied in this paper does not consider the decision variables associated with controllable power system devices such as high-voltage DC (HVDC) lines, tap-changing and phase-shifting transformers, and switched-shunt devices. Incorporating the continuous and potentially discrete decision variables necessary for modeling these devices in convex relaxations of the OPF problem is an area of ongoing research. For instance, a convex, second-order cone programming formulation for HVDC lines is available in~\cite{socp_hvdc}. One possibility is modeling discrete variables as polynomial equality constraints (e.g., $\phi \in \left\lbrace 0,1 \right\rbrace$ is equivalent to \mbox{$\phi^2 - \phi = 0$}), which can be incorporated in the moment relaxations discussed in this paper. This is a promising direction for future work.

\section{Moment Relaxations}
\label{l:msdp_overview}


\subsection{Overview}

The OPF problem~\eqref{opf} is comprised of polynomial functions of the voltage components $V_d$ and $V_q$ and can therefore be solved using moment relaxations~\cite{lasserre2001,lasserre_book}. We next present moment relaxations of the OPF problem~\eqref{opf}. The material in this section builds on~\cite{pscc2014}. More detailed descriptions of moment relaxations are available in~\cite{lasserre2001} and~\cite{lasserre_book}, and application of moment relaxations to the OPF problem is independently proposed in~\cite{patrick} and~\cite{ibm_paper}.

Polynomial optimization problems, such as the OPF problem, are a special case of \emph{generalized moment problems}~\cite{lasserre_book}. Global solutions to generalized moment problems can be approximated using moment relaxations that are formulated as semidefinite programs. For polynomial optimization problems with bounded variables, such as OPF problems, the approximation approaches the global solution(s) as the relaxation order increases~\cite{lasserre_book}. While moment relaxations can find all global solutions to polynomial optimization problems, we focus on problems with a single global optimum.

Formulating the moment relaxations requires several definitions. Define the vector $\hat{x} = \begin{bmatrix} V_{d1} & V_{d2} & \ldots & V_{qn} \end{bmatrix}^\intercal$, which contains all first-order monomials of the decision variables in~\eqref{opf}. Given a vector \mbox{$\alpha \in \mathbb{N}^{2n}$} representing monomial exponents, the expression $\hat{x}^\alpha = V_{d1}^{\alpha_1}V_{d2}^{\alpha_2}\cdots V_{qn}^{\alpha_{2n}}$ defines the monomial associated with $\hat{x}$ and $\alpha$. A polynomial $g\left(\hat{x}\right)$ can be expressed as

\begin{equation}\label{gpoly}
g\left(\hat{x}\right) \triangleq \sum_{\alpha \in \mathbb{N}^{2n}} g_{\alpha} \hat{x}^{\alpha}
\end{equation}

\noindent where $g_{\alpha}$ is the scalar coefficient corresponding to the monomial $\hat{x}^{\alpha}$.

Next define a linear functional $L_y\left\lbrace g\right\rbrace$:

\begin{equation}\label{L}
L_y\left\lbrace g\right\rbrace \triangleq \sum_{\alpha \in \mathbb{N}^{2n}} g_{\alpha} y_{\alpha}
\end{equation}

\noindent This functional replaces the monomials $\hat{x}^\alpha$ in a polynomial function $g\left(\hat{x}\right)$ with scalar variables $y_{\alpha}$. When $g\left(\hat{x}\right)$ is a matrix, the functional $L_y\left\lbrace g \right\rbrace$ is applied to each element of $g\left(\hat{x}\right)$.

Consider, for example, the vector $\hat{x} = \begin{bmatrix}V_{d1} & V_{d2} & V_{q2} \end{bmatrix}^\intercal$ corresponding to the voltage components of a two-bus system, where the angle reference constraint~\eqref{opf_Vref} is used to eliminate $V_{q1}$, and the polynomial $g\left(\hat{x}\right) = -\left(0.95\right)^2 + f_{V2}\left(V_d,V_q \right) = -\left(0.95\right)^2 + V_{d2}^2 + V_{q2}^2$. (The constraint $g\left(\hat{x}\right) \geq 0$ forces the voltage magnitude at bus~2 to be greater than or equal to 0.95~per unit.) Then $L_y\left\lbrace g\right\rbrace = -\left(0.95\right)^2y_{000} + y_{020} + y_{002}$. Thus, $L\left\lbrace g \right\rbrace$ converts a polynomial $g\left(\hat{x}\right)$ to a linear function of $y$.

The order-$\gamma$ relaxation forms a vector $x_\gamma$ composed of all monomials of the voltage components up to order $\gamma$:

\begin{align} \nonumber
x_\gamma \triangleq & \left[ \begin{array}{ccccccc} 1 & V_{d1} & \ldots & V_{qn} & V_{d1}^2 & V_{d1}V_{d2} & \ldots \end{array} \right. \\ \label{x_d}
& \qquad \left.\begin{array}{cccccc} \ldots & V_{qn}^2 & V_{d1}^3 & V_{d1}^2 V_{d2} & \ldots & V_{qn}^\gamma \end{array}\right]^\intercal
\end{align}

We now define \emph{moment} and \emph{localizing} matrices. The symmetric moment matrix $\mathbf{M}_\gamma\left(y\right)$ has entries $y_\alpha$ corresponding to all monomials $\hat{x}^\alpha$ up to order $2\gamma$:

\begin{equation}\label{momentmat}
\mathbf{M}_\gamma \left\lbrace y \right\rbrace \triangleq L_y\left\lbrace x_\gamma^{\vphantom{\intercal}} x_\gamma^\intercal\right\rbrace
\end{equation}



Symmetric localizing matrices\footnote{The terminology \emph{localizing matrix} is adopted from~\cite{lasserre2001,lasserre_book}.} are defined for each constraint of~\eqref{opf}. The localizing matrices consist of linear combinations of the moment matrix entries $y$. Each polynomial constraint of the form $f\left(\hat{x}\right) - a \geq 0$ in~\eqref{opf} (e.g., $f_{V2}\left(\hat{x}\right) - V_2^{\min} \geq 0$) corresponds to the localizing matrix

\begin{equation}\label{localizing} \small
\mathbf{M}_{\gamma-\beta}\left\lbrace \left(f\left(\hat{x}\right) - a\right) y \right\rbrace \triangleq L_y\left\lbrace\left(f\left(\hat{x}\right) - a\right) x_{\gamma-\beta}^{\vphantom{\intercal}} x_{\gamma-\beta}^\intercal \right\rbrace
\end{equation}

\noindent where the polynomial $f$ has degree $2\beta$. Example moment and localizing matrices for the second-order relaxation of a two-bus system are presented in~\eqref{2busMomentMat} and~\eqref{2busLocalizingMat}, respectively.

The order-$\gamma$ moment relaxation of~\eqref{opf} is

\begin{subequations}\small
\label{msdp_opf}
\begin{align}
\label{msdp_obj}& \min_{y} L_y\left\lbrace \sum_{k \in \mathcal{G}} f_{Ck} \right\rbrace \qquad \mathrm{subject\; to} \hspace{-150pt} &  \\
\label{msdp_Pmin} & \quad \mathbf{M}_{\gamma-1}\left\lbrace \left(f_{Pk} - P_k^{\min}\right) y \right\rbrace \succeq 0 & \forall k\in\mathcal{N}\\
\label{msdp_Pmax} & \quad \mathbf{M}_{\gamma-1}\left\lbrace \left(P_k^{\max} - f_{Pk} \vphantom{P_k^{\min}}\right) y \right\rbrace \succeq 0 & \forall k\in\mathcal{N}\\
\label{msdp_Qmin} & \quad \mathbf{M}_{\gamma-1}\left\lbrace \left(f_{Qk} - Q_k^{\min}\right) y \right\rbrace \succeq 0 & \forall k\in\mathcal{N}\\
\label{msdp_Qmax} & \quad \mathbf{M}_{\gamma-1}\left\lbrace \left(Q_k^{\max} - f_{Qk}  \vphantom{P_k^{\min}}\right) y \right\rbrace \succeq 0 & \forall k\in\mathcal{N}\\
\label{msdp_Vmin} & \quad \mathbf{M}_{\gamma-1}\left\lbrace \left(f_{Vk} - V_k^{\min}\right) y \right\rbrace \succeq 0 & \forall k\in\mathcal{N}\\
\label{msdp_Vmax} & \quad \mathbf{M}_{\gamma-1}\left\lbrace \left(V_k^{\max} - f_{Vk}  \vphantom{P_k^{\min}}\right) y \right\rbrace \succeq 0 & \forall k\in\mathcal{N} \\
\label{msdp_Slm} & \quad \mathbf{M}_{\gamma-2}\left\lbrace \left(S_{lm}^{\max} - f_{Slm} \vphantom{P_k^{\min}}\right) y \right\rbrace \succeq 0 & \forall \left(l,m\right)\in\mathcal{L} \\
\label{msdp_Sml} & \quad \mathbf{M}_{\gamma-2}\left\lbrace \left(S_{lm}^{\max} - f_{Sml} \vphantom{P_k^{\min}}\right) y \right\rbrace \succeq 0 & \forall \left(l,m\right)\in\mathcal{L} \\
\label{msdp_Msdp} & \quad \mathbf{M}_\gamma \left(y\right) \succeq 0 & \\
\label{msdp_y0} & \quad y_{00\ldots 0} = 1 & \\
\label{msdp_Vref} & \quad y_{0\ldots00\eta0\ldots0} = 0 & \eta = 1,\ldots,2\gamma
\end{align}
\end{subequations}

\begin{figure}[b]
\centering
\includegraphics[totalheight=0.09\textheight]{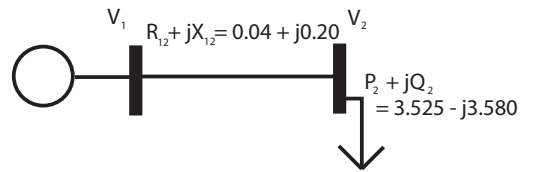}
\caption{Two-Bus System from ~\cite{bukhsh2011}}
\label{f:twobussystem}
\end{figure}

\begin{figure*}[t]
\begin{equation}\label{2busMomentX} 
x_2 = \left[\begin{array}{cccccccccc} 1 & V_{d1} & V_{d2} & V_{q2} & V_{d1}^2 & V_{d1}V_{d2} & V_{d1}V_{q2} & V_{d2}^2 & V_{d2}V_{q2} & V_{q2}^2\end{array}\right]^\intercal
\end{equation}
\vspace{3pt}
\begin{equation}\label{2busMomentMat}
\mathbf{M}_2 \left(y \right) = L_y\left(x_2 x_2^\intercal\right) = \left[\begin{array}{c|ccc|cccccc} 
y_{000} & y_{100} & y_{010} & y_{001} & y_{200} & y_{110} & y_{101} & y_{020} & y_{011} & y_{002} \\\hline
y_{100} & y_{200} & y_{110} & y_{101} & y_{300} & y_{210} & y_{201} & y_{120} & y_{111} & y_{102} \\
y_{010} & y_{110} & y_{020} & y_{011} & y_{210} & y_{120} & y_{111} & y_{030} & y_{021} & y_{012} \\
y_{001} & y_{101} & y_{011} & y_{002} & y_{201} & y_{111} & y_{102} & y_{021} & y_{012} & y_{003} \\ \hline
y_{200} & y_{300} & y_{210} & y_{201} & y_{400} & y_{310} & y_{301} & y_{220} & y_{211} & y_{202} \\ 
y_{110} & y_{210} & y_{120} & y_{111} & y_{310} & y_{220} & y_{211} & y_{130} & y_{121} & y_{112} \\
y_{101} & y_{201} & y_{111} & y_{102} & y_{301} & y_{211} & y_{202} & y_{121} & y_{112} & y_{103} \\
y_{020} & y_{120} & y_{030} & y_{021} & y_{220} & y_{130} & y_{121} & y_{040} & y_{031} & y_{022} \\
y_{011} & y_{111} & y_{021} & y_{012} & y_{211} & y_{121} & y_{112} & y_{031} & y_{022} & y_{013} \\ 
y_{002} & y_{102} & y_{012} & y_{003} & y_{202} & y_{112} & y_{103} & y_{022} & y_{013} & y_{004} 
\end{array}\right]
\end{equation}
\begin{align}\label{2busLocalizingMat}  \nonumber
& \mathbf{M}_{1}\left\lbrace\left(f_{V2} - \left(0.95\right)^2\right) y \right\rbrace = \\
& \quad \left[\begin{array}{c|ccc}
y_{020} + y_{002} - \left(0.95\right)^2y_{000} & y_{120} + y_{102} - \left(0.95\right)^2y_{100} & y_{030} + y_{012} - \left(0.95\right)^2y_{010} & y_{021} + y_{003} - \left(0.95\right)^2y_{001}  \\ \hline
y_{120} + y_{102} - \left(0.95\right)^2y_{100} & y_{220} + y_{202} - \left(0.95\right)^2y_{200} & y_{130} + y_{112} - \left(0.95\right)^2y_{110} & y_{121} + y_{103} - \left(0.95\right)^2y_{101} \\
y_{030} + y_{012} - \left(0.95\right)^2y_{010} & y_{130} + y_{112} - \left(0.95\right)^2y_{110} & y_{040} + y_{022} - \left(0.95\right)^2y_{020} & y_{031} + y_{013} - \left(0.95\right)^2y_{011} \\
y_{021} + y_{003} - \left(0.95\right)^2y_{001} & y_{121} + y_{103} - \left(0.95\right)^2y_{101} & y_{031} + y_{013} - \left(0.95\right)^2y_{011} & y_{022} + y_{004} - \left(0.95\right)^2y_{002}  \\
\end{array}\right]
\end{align}
\vspace*{3pt}
\hrule
\vspace*{0pt}
\end{figure*}

\noindent where $\succeq 0$ indicates that the corresponding matrix is positive semidefinite. The moment relaxation is thus a semidefinite program. (A dual form of the moment relaxation is a sum-of-squares program~\cite{lasserre_book}.) Note that the constraint~\eqref{msdp_y0} enforces the fact that $x^{0} = 1$. The constraint~\eqref{msdp_Vref} corresponds to the angle reference constraint~\eqref{opf_Vref}; the $\eta$ in \eqref{msdp_Vref} is in the index $n+1$, which corresponds to the variable $V_{q1}$. Note that the angle reference can alternatively be used to eliminate all terms corresponding to $V_{q1}$ to reduce the size of the semidefinite program.

\subsection{Two-Bus Example}
\label{l:twobusexample}

We next present an illustrative two-bus example problem from~\cite{bukhsh2011}. Fig.~\ref{f:twobussystem} gives the system's one-line diagram assuming a 100~MVA base power. The generator at bus~1 has no limits on active or reactive outputs and there is no line-flow limit. Bus~1 voltage magnitude is in the range $\left[0.95, 1.05\right]$~per unit, while bus~2 voltage magnitude is in the range $\left[0.95,1.02\right]$~per unit. Specify a \$1/MWh cost of active power generation at bus~1. With three degrees of freedom, the entire feasible space of the two-bus OPF problem can be visualized in three dimensions as shown in Section~IV of~\cite{pscc2014}.

Note that while this system has a radial network topology, the OPF problem does not satisfy the sufficient conditions for exactness of the semidefinite relaxation described in~\cite{low_tutorial1,low_tutorial2}, and the first-order moment relaxation is, in fact, not exact for this problem.

The second-order relaxation has the vector $x_2$ given in~\eqref{2busMomentX} and the moment matrix in~\eqref{2busMomentMat}. The localizing matrix corresponding to the voltage magnitude constraint $f_{V2}\left(V_d,V_q \right) \geq \left(0.95\right)^2$ in~\eqref{opf_V} is given in~\eqref{2busLocalizingMat}. Note that the angle reference constraint~\eqref{opf_Vref} is used to eliminate $V_{q1}$ so that $\hat{x} = \begin{bmatrix}V_{d1} & V_{d2} & V_{q2} \end{bmatrix}^\intercal$.

The second-order moment relaxation yields the global solution of $V = \begin{bmatrix} 0.950 & 0.416 - j0.893 \end{bmatrix}^\intercal$~per unit. This corresponds to active and reactive power generation at bus~1 of 456.6 MW and 162.3 MVAr, respectively, and an operating cost of 456.55~\$/MWh.

\subsection{Implementation Details}

This section next discusses several implementation details for the moment relaxations, including removing linear dependency resulting from equality-constrained polynomials, recovering the global solution(s) to the OPF problem from a solution to the moment relaxation, elimination of unnecessary portions of the moment relaxation's matrices, and reformulation of the cost function and apparent-power line-flow constraints. Note that, in contrast to the example in Section~\ref{l:twobusexample}, the notation in this section does not use the angle reference constraint $V_{q1}=0$ to eliminate $V_{q1}$.

\subsubsection{Equality-Constrained Polynomials}

Given that the localizing matrices are symmetric, the constraints $\mathbf{M}_{\gamma-1}\left\lbrace \left(f\left(\hat{x}\right) - a\right) y \right\rbrace \succeq 0$ and $\mathbf{M}_{\gamma-1}\left\lbrace \left(a - f\left(\hat{x}\right)\right) y \right\rbrace \succeq 0$ imply that the matrix $\mathbf{M}_{\gamma-1}\left\lbrace \left(f\left(\hat{x}\right) - a\right) y \right\rbrace = 0$. Thus, all entries of a localizing matrix corresponding to an equality-constrained polynomial $f\left(\hat{x}\right) = a$ (e.g., the power flow constraints at load buses) are zero. That is, equality-constrained polynomials result in equality constraints that are linear in the variables $y$ rather than positive semidefinite matrix constraints.

The non-uniqueness of the entries of the localizing matrices (see, e.g., the repeated terms in~\eqref{2busLocalizingMat} due to the matrix's symmetry) creates linearly dependent equality constraints which can introduce numerical difficulties. To eliminate these redundant constraints, the localizing matrix constraints for equality-constrained polynomials are therefore replaced by the equivalent vector of constraints $L_y\left\lbrace g\left(\hat{x} \right) x_{\gamma-1} \right\rbrace = 0$, where $g\left(\hat{x}\right)$ denotes any equality-constrained polynomial in~\eqref{opf}.

\subsubsection{Solution Extraction}
The \mbox{order-$\gamma$} moment relaxation yields a single global solution if $\mathrm{rank}\left(\mathbf{M}_{\gamma}\left(y\right)\right) = 1$. The global solution $x^\ast$ to the OPF problem~\eqref{opf} is then determined by a spectral decomposition of the diagonal block of the moment matrix corresponding to the second-order terms. Specifically, let $\eta$ be a unit-length eigenvector corresponding to the non-zero eigenvalue $\lambda$ from the diagonal block of the moment matrix corresponding to the second-order monomials (i.e., $\left[\mathbf{M}_2\right]_{\left(2:k,2:k\right)}$, where $k = 2n+1$ and subscripts indicate the vector entries in MATLAB notation). Then the vector $V^\ast = \sqrt{\lambda} \left(\eta_{1:n} + j \eta_{\left(n+1\right):2n}\right)$ is the globally optimal voltage phasor vector.

If $\mathrm{rank}\left(\mathbf{M}_{\gamma}\left(y\right) \right) > 1$, there are either multiple global solutions (i.e., there are multiple points in the feasible space of the original non-convex OPF problem with the same globally optimal objective value) requiring the solution extraction procedure described in Section 5.3.1 of~\cite{lasserre_book}\footnote{If $\mathrm{rank}\left(\mathbf{M}_{\gamma - 1}\left(y \right) \right) = \mathrm{rank}\left(\mathbf{M}_{\gamma} \left(y\right)\right)$, then there are at least $\mathrm{rank}\left(\mathbf{M}_{\gamma} \left(y\right)\right)$ globally optimal solutions to the OPF problem. The globally optimal decision variable vectors can be extracted using Algorithm~4.2 in~\cite{lasserre_book}, which only uses linear algebra operations. \\ \indent While practical OPF problems can have multiple \emph{local} solutions, we expect that multiple \emph{global} solutions are uncommon. Reference~\cite{lavaei_tps} uses example OPF problems with multiple global solutions to show that the OPF problem is, in general, NP-hard. As discussed in~\cite{pscc2014}, these example problems are atypical, but are useful for exploring the limits of the moment relaxation approach.} or the order-$\gamma$ moment relaxation is not exact and only yields a lower bound on the objective value. If the order-$\gamma$ moment relaxation is not exact, the order-$\left(\gamma + 1\right)$ moment relaxation will improve the lower bound and may give a global solution.


\subsubsection{Elimination of Unnecessary Terms}

Since the polynomials in~\eqref{opf} are composed solely of constant, second-, and fourth-order monomials, off-diagonal blocks of the moment matrix corresponding to odd-order monomials are not required. Further, all terms in the off-diagonal blocks of the moment matrix corresponding to even-order monomials are duplicated in the diagonal blocks and are therefore unnecessary. (See~\eqref{2busMomentMat} for an illustration of this matrix partitioning.) Thus, positive semidefinite matrix constraints are only enforced for the diagonal blocks of the moment matrix corresponding to even-order monomials (i.e., $y_{\alpha}$ such that $\sum_{i=1}^{2n} \alpha_i$ is even). Similarly, positive semidefinite constraints are only applied to the diagonal blocks of the localizing matrices which correspond to the even-order monomials of the matrix $x_{\gamma-\beta}^{\vphantom{\intercal}} x_{\gamma-\beta}^\intercal$ (e.g., the diagonal blocks of~\eqref{2busLocalizingMat}). By reducing the size of the semidefinite program, this decreases the relaxation's computational burden. 

After elimination of terms corresponding to the first-order monomials, the moment matrix for the first-order relaxation contains only the second-order monomials. Further, the localizing ``matrices'' are in fact positivity constraints on scalars.\footnote{For the first-order moment relaxation, the localizing matrix for the polynomial constraint $f\left(x\right) \geq 0$ is $L_y\left( f\left(x\right)x_0^{{{\vphantom{\intercal}}}}x_0^\intercal\right) = L_y\left(f\left(x\right)\cdot 1\cdot 1^\intercal\right) = L_y\left(f\left(x\right)\right) \geq 0$ (i.e., a positive scalar constraint).} Thus, the first-order relaxation is equivalent to the semidefinite relaxation of~\cite{lavaei_tps}. The higher-order moment relaxations are generalizations of the semidefinite relaxation of~\cite{lavaei_tps}.

\subsubsection{Quadratic Cost Function and Apparent-Power Line-FLow Limits}

The order $\gamma$ of the moment relaxation must be greater than or equal to half of the degree of any polynomial in the OPF problem~\eqref{opf}. Relaxations of all polynomials can then be written as linear functions of the entries of $\mathbf{M}_\gamma$. For instance, the OPF problem with a linear cost function and without apparent-power line-flow limits requires $\gamma \geq 1$. Although direct implementation of~\eqref{opf} requires $\gamma \geq 2$ due to the fourth-order polynomials in the cost function~\eqref{opf_obj} and apparent-power line-flow limits~\eqref{opf_Slm}-\eqref{opf_Sml}, these fourth-order polynomials can be rewritten as second-order polynomials using a Schur complement formulation~\cite{lavaei_tps}. Specifically, rather than~\eqref{msdp_Slm}-\eqref{msdp_Sml}, enforce the constraints



\begin{subequations}
\label{lineflow_schur}
\begin{align}
 & 
\begin{bmatrix}
-\left(S_{lm}^{\mathrm{max}}\right)^2 & L_y\left\lbrace f_{Plm} \right\rbrace  & L_y\left\lbrace f_{Qlm} \right\rbrace \\
L_y\left\lbrace f_{Plm} \right\rbrace & -1 & 0 \\
L_y\left\lbrace f_{Qlm} \right\rbrace & 0 & -1
\end{bmatrix} \preceq 0 &  & \forall \left(l,m\right) \in \mathcal{L}
\end{align}
\begin{align}
\begin{bmatrix}
-\left(S_{lm}^{\mathrm{max}}\right)^2 & L_y\left\lbrace f_{Pml} \right\rbrace  & L_y\left\lbrace f_{Qml} \right\rbrace \\
L_y\left\lbrace f_{Pml} \right\rbrace & -1 & 0 \\
L_y\left\lbrace f_{Qml} \right\rbrace & 0 & -1
\end{bmatrix} \preceq 0 &  & \forall \left(l,m\right) \in \mathcal{L}
\end{align}
\end{subequations} 


Similarly, define new variables $\alpha_k$ for each generator $k\in\mathcal{G}$ and replace the quadratic cost function in~\eqref{msdp_obj} with $\sum_{k\in\mathcal{G}} \alpha_k$ and the additional constraint




\begin{align}
\label{quad_schur} & \begin{bmatrix}
c_{k1} L_y\left\lbrace f_{Pk}\right\rbrace + c_{k0} - \alpha_k & \sqrt{c_{k2}} L_y\left\lbrace f_{Pk}\right\rbrace \\
\sqrt{c_{k2}} L_y\left\lbrace f_{Pk}\right\rbrace & -1
\end{bmatrix} \preceq 0 & & \forall k \in \mathcal{G}
\end{align}

\noindent Note that a second-order cone programming (SOCP) formulation can also be employed to represent the quadratic cost function and apparent-power line-flow limits~\cite{andersen2014}.

The OPF problem reformulated using~\eqref{lineflow_schur} and~\eqref{quad_schur} only requires $\gamma \geq 1$. Use of both the Schur complement formulations and direct implementation for the apparent-power line-flow constraints~\eqref{opf_Slm}-\eqref{opf_Sml} and quadratic cost function~\eqref{msdp_obj} generally gives superior results for $\gamma \geq 2$ as compared to implementing either the Schur complement or direct formulations separately. That is, when possible, enforce both \eqref{msdp_Slm}-\eqref{msdp_Sml} and \eqref{lineflow_schur} for the apparent-power line-flow constraints and include the constraints~\eqref{quad_schur} and \mbox{$\alpha_k = L\left\lbrace f_{Ck}\right\rbrace \; \forall k\in \mathcal{G}$} while minimizing $\sum_{k \in \mathcal{G}}\alpha_k$ for the quadratic cost function).

\section{Exploiting Sparsity in Moment Relaxations}
\label{l:msdp_sparse}


The moment relaxations globally solve a broader class of OPF problems than existing convex relaxations~\cite{pscc2014,patrick,ibm_paper}. However, the superior capabilities of the moment relaxations have a computational cost: the semidefinite program needed to solve the moment relaxation quickly becomes computationally intractable with both increasing problem size and relaxation order. Direct implementation of the formulation presented in Section~\ref{l:msdp_overview} is computationally tractable for a second-order relaxation of OPF problems with up to ten buses. 


Solving larger OPF problems requires exploiting power system sparsity. Similar to methods for existing semidefinite relaxations~\cite{jabr11,molzahn_holzer_lesieutre_demarco-large_scale_sdp_opf}, matrix completion decomposition is applicable to moment relaxations of polynomial optimization problems~\cite{waki2006}. This decomposition extends computational tractability to OPF problems with approximately forty buses.

Solving OPF problems with more than forty buses requires exploiting the observation that the first-order relaxation is sufficient for large regions of typical OPF problems~\cite{hicss2014}. By selectively applying second- and third-order relaxations to specific buses or small groups of buses, larger OPF problems become computationally tractable. With standard semidefinite programming solvers, OPF problems of up to 300 buses can be solved.

After reviewing the matrix completion decomposition~\cite{waki2006}, this section  proposes a method for selectively applying the moment relaxation and presents a heuristic method for determining where to apply higher-order relaxations.

\subsection{Matrix Completion Decomposition}
\label{l:matrix_completion}

The matrix completion decomposition, which is adopted from~\cite{waki2006}, exploits power system sparsity. The decomposition relies on a matrix completion theorem~\cite{gron1984,fukuda2001} which draws on graph theory. Several graph theoretic definitions are necessary for understanding the matrix completion theorem. A \emph{clique} is a subset of the graph nodes for which each node in the clique is connected to all other nodes in the clique. A \emph{maximal clique} is a clique that is not a proper subset of another clique. Denote the set of maximal cliques by $\mathcal{M}$, with $\mathcal{M}_m$ representing the set of buses associated with the $m^{th}$ maximal clique. A graph is \emph{chordal} if each cycle of length four or more nodes has a chord, which is an edge connecting two nodes that are not adjacent in the cycle.

The graph in question for the moment relaxations of the OPF problem is defined with a set of nodes $\mathcal{\hat{N}}$ and a set of undirected edges $\mathcal{\hat{L}}$. This graph is derived from the power system network. The set of nodes is equal to the set of buses in the power system (i.e., $\mathcal{\hat{N}} = \mathcal{N} = \left\lbrace 1,\ldots,n\right\rbrace$). The set of edges $\mathcal{\hat{L}}$ is a superset of the topology of the power system network $\mathcal{L}$. Define $\mathcal{\bar{N}}_k$ as the subset of buses connected to bus~$k$ in the power system network (i.e., $\mathcal{\bar{N}}_k = \left\lbrace i \,\left|\, \left(i,k \right) \in \mathcal{L} \right. \right\rbrace$). For each bus~$k$, add to $\mathcal{\hat{L}}$ all edges between each bus in $\mathcal{\bar{N}}_k$. That is, all neighboring buses of each bus are connected in $\mathcal{\hat{L}}$.\footnote{This is a subtle but important difference from the matrix completion decompositions in~\cite{jabr11} and~\cite{molzahn_holzer_lesieutre_demarco-large_scale_sdp_opf} that directly use the graph with nodes $\mathcal{N}$ and edges $\mathcal{L}$ from the power system network. Matrix completion decompositions for the higher-order moment formulations use the graph defined by $\hat{\mathcal{N}}$ and $\hat{\mathcal{L}}$ so that each bus belongs to at least one maximal clique that also contains each of that bus' neighbors, and thus all variables necessary for the higher-order moment constraints are well-defined~\cite{waki2006}.}

The maximal cliques of a chordal graph can be determined in linear time~\cite{tarjan}. However, identifying the maximal cliques of a non-chordal graph is an NP-hard problem. Since realistic power networks are generally not chordal, we use a chordal extension technique which adds edges to $\mathcal{\hat{L}}$ to obtain a chordal super-graph denoted as $\mathcal{\hat{L}}_{ch}$. To form the chordal extension, denote as $\mathbf{D}$ the adjacency matrix of the graph defined by $\mathcal{\hat{N}}$ and $\mathcal{\hat{L}}$. The chordal extension is then determined using a Cholesky factorization of $\mathbf{D} + \mathbf{I}$, where $\mathbf{I}$ is an $n\times n$ identity matrix. The off-diagonal sparsity pattern of $\mathrm{chol}\left(\mathbf{D}+\mathbf{I}\right)$ provides a chordal extension $\mathcal{\hat{L}}_{ch}$. An approximate minimum-degree permutation of the buses~\cite{davis2004} is employed to reduce the number of added edges in $\mathcal{\hat{L}}_{ch}$ relative to $\mathcal{\hat{L}}$.

The matrix completion theorem can now be stated. Let $\mathbf{W}$ be a symmetric matrix with partial information (i.e., not all entries of $\mathbf{W}$ have known values) with an associated undirected graph. (For the moment relaxation, the graph in question has nodes $\mathcal{\hat{N}}$ and edges $\mathcal{\hat{L}}_{ch}$.) The matrix $\mathbf{W}$ can be completed to a positive semidefinite matrix (i.e., the unknown entries of $\mathbf{W}$ can be chosen such that $\mathbf{W} \succeq 0$) if and only if the submatrices associated with each of the maximal cliques of the graph defined by $\mathbf{W}$ are all positive semidefinite.

The matrix completion theorem allows replacing the single large positive semidefinite constraint on the moment matrix~\eqref{msdp_Msdp} with constraints on many smaller matrices:

\begin{align}\label{moment_decomposed}
\mathbf{M}_\gamma^{\mathcal{M}_m} \left( y \right) \triangleq L_y\left\lbrace \left(x_\gamma^{\mathcal{M}_m}\right) \left(x_\gamma^{\mathcal{M}_m}\right)^\intercal \right\rbrace \succeq 0 & & m = 1,\ldots,\left|\mathcal{M}\right|
\end{align}

\noindent where $x_\gamma^{\mathcal{M}_m}$ is the subset of $x_\gamma$ corresponding to the buses in $\mathcal{M}_m$ and $\left|\mathcal{M}\right|$ is the number of maximal cliques in the graph. 

Similarly, the localizing matrices in \eqref{msdp_Pmin}-\eqref{msdp_Sml} are each replaced by a single smaller matrix. Each bus~$k$ is associated with a single \emph{smallest covering} maximal clique $m_k \in \left\lbrace 1,\ldots,\left|\mathcal{M}\right|\right\rbrace$ (i.e., the maximal clique with least number of buses that completely contains bus $k$ and its neighbors in the power system network). By construction of $\mathcal{\hat{L}}$, each bus and its neighbors will be entirely contained in at least one maximal clique. Form the localizing matrix

\begin{align}\nonumber
& \mathbf{M}_{\gamma-\beta}^{\mathcal{M}_{m_k}}\left\lbrace \left(f\left(\hat{x} \right) - a\right) y \right\rbrace \\ \label{localizing_decomposed}
& \qquad \triangleq L_y\left\lbrace\left(f\left(\hat{x} \right) - a\right) \left(x_{\gamma-\beta}^{\mathcal{M}_{m_k}}\right) \left(x_{\gamma-\beta}^{\mathcal{M}_{m_k}}\right)^\intercal \right\rbrace
\end{align}

\noindent where $f\left(\hat{x} \right) - a \geq 0$ denotes a generic polynomial constraint in~\eqref{opf} with order $2\beta$ associated with bus~$k$.

Since the maximal cliques have non-empty intersection (i.e., contain some of the same buses), different decomposed moment matrices may contain elements that refer to a common element in the original moment matrix. The decomposed optimization problem must be formulated such that these shared elements are equal.

The solution to the decomposed formulation consists of many matrices. The globally optimal voltage vector solution to the OPF problem can be recovered if each moment matrix satisfies a rank condition. Specifically, the diagonal blocks corresponding to the second-order monomials in all moment matrices must have rank one. If the rank condition is satisfied, the method described in~\cite{molzahn_holzer_lesieutre_demarco-large_scale_sdp_opf} may be used to recover the globally optimal voltage vector. 

If any of the decomposed moment matrices does not satisfy the rank condition, the decomposed moment relaxation does not yield a solution to the OPF problem. Failure to satisfy the rank condition may either indicate that the moment relaxation is not exact or that there are multiple global solutions. In the former case, the objective value from the moment relaxation serves as a lower bound on the objective value of the OPF problem~\eqref{opf} and increasing the relaxation order may result in a relaxation that yields a global solution.\footnote{The lower bound may in fact be the global minimum objective value without the relaxation providing globally optimal decision variables (i.e., a ``hidden'' rank one solution~\cite{bitar_allerton2013}).} In the latter case, unlike the formulation in Section~\ref{l:msdp_overview} which can recover multiple global solutions using the method described in~\cite{lasserre_book}, there is only limited ability to recover multiple global solutions to the decomposed moment relaxation. Adding a small perturbation to the objective function may result in recovery of a single global solution~\cite{waki2006}.

\subsection{Selective Application of Higher-Order Constraints}
\label{l:selective_application}

The matrix completion decomposition described in Section~\ref{l:matrix_completion} significantly reduces the size of the semidefinite program for large, sparse power networks. With this decomposition, second-order relaxations of OPF problems with up to approximately forty buses are computationally tractable.

Solving larger OPF problems is accomplished by exploiting the observation that the first-order relaxation is sufficient for large regions of typical OPF problems. A voltage vector is obtained from the closest matrix that satisfies the rank condition (i.e., the rank-one matrix with smallest Frobenius-norm difference to the moment matrix from the relaxation's solution). This matrix is determined using an eigen decomposition of the higher-rank diagonal block of the moment matrix corresponding to the second-order monomials. A power injection ``mismatch'' is determined by comparing the value of the power injections calculated from the higher-rank matrix ($L_y\left\lbrace f_{Pk} \right\rbrace$ and $L_y\left\lbrace f_{Qk} \right\rbrace$) to the power injections implied by the voltage vector from the closest matrix satisfying the rank condition.\footnote{At load buses, the power injection mismatches are equal to the difference between the specified load demands and the power injections implied by the closest rank one matrix. At generator buses, the mismatches are equal to the difference between the power injections derived from the localizing matrices (i.e., the elements in the (1,1) position of~\eqref{msdp_Pmin} and~\eqref{msdp_Qmin} plus the load demands) and the power injections implied by the closest rank one matrix.}

\begin{figure}[t]
\centering
\includegraphics[totalheight=0.26\textheight]{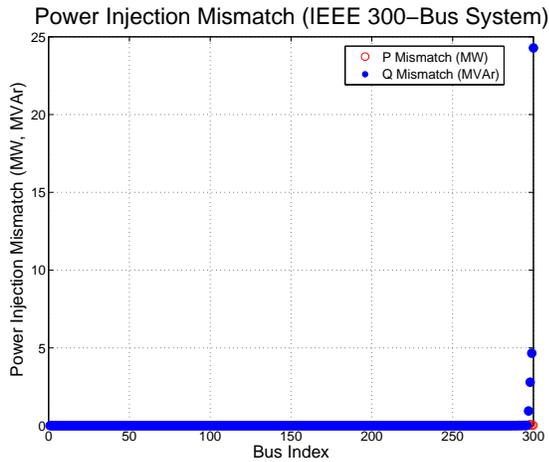}
\caption{Power Injection Mismatches for the First-Order Relaxation of the IEEE 300-Bus System~\cite{hicss2014}}\label{f:case300mismatch} 
\end{figure}

First-order relaxations typically yield voltage vectors that have small power injection mismatches at the majority of buses while a few buses have large mismatch~\cite{hicss2014}. For example, Fig.~\ref{f:case300mismatch} shows the power injection mismatches, sorted in increasing order, resulting from the first-order moment relaxation of the IEEE 300-bus system. This suggests that selective application of higher-order relaxations to specified buses may be sufficient to globally solve larger OPF problems. 

To selectively apply the higher-order constraints, each bus~$k$ has an associated relaxation order $\gamma_k$ rather than a uniform order for the entire OPF problem. (A heuristic for specifying the $\gamma_k$ values is presented in Section~\ref{l:heuristic}.) Determine the relaxation order for each maximal clique $\mu = 1,\ldots,\left|\mathcal{M}\right|$, denoted as $\hat{\gamma}_\mu$, which is the highest relaxation order of any bus for which $\mu$ is the smallest covering maximal clique.\footnote{Note that $\hat{\gamma}_\mu$ is not necessarily the highest relaxation order of the buses in maximal clique $\mu$. The maximal clique $\mu$ may contain buses which are also contained in other maximal cliques. The order $\hat{\gamma}_\mu$ is determined by the buses for which $\mu$ is the \emph{smallest covering} maximal clique (i.e., $\mu$ is the smallest clique to contain that bus and all of its neighbors).} The decomposed moment matrix constraints in~\eqref{moment_decomposed} are formed according to this order:

\begin{align}\label{moment_decomposed_selective}
\mathbf{M}_{\hat{\gamma}_\mu}^{\mathcal{M}_\mu} \left( y \right) \triangleq L_y\left\lbrace \left(x_{\hat{\gamma}_\mu}^{\mathcal{M}_\mu}\right) \left(x_{\hat{\gamma}_\mu}^{\mathcal{M}_\mu}\right)^\intercal \right\rbrace \succeq 0 & & \mu = 1,\ldots,\left|\mathcal{M}\right|
\end{align}

The localizing matrices for the power injection and voltage magnitude constraints \eqref{msdp_Pmin}-\eqref{msdp_Vmax} at bus~$k$ are constructed according to the corresponding bus order $\gamma_k$, while the localizing matrices for the apparent-power line-flow constraints are constructed according to the highest order of either terminal bus. That is, for each constraint $f\left(\hat{x} \right) - a \geq 0$ with order $2\beta$ in \eqref{opf}, create the localizing matrix constraint

\begin{align}\nonumber
& \mathbf{M}_{\gamma_k-\beta}^{\mathcal{M}_m}\left\lbrace \left(f\left(\hat{x} \right) - a\right) y \right\rbrace \\ \label{localizing_decomposed_selective}
& \qquad \triangleq L_y\left\lbrace\left(f\left(\hat{x} \right) - a\right) \left(x_{\gamma_k-\beta}^{\mathcal{M}_m}\right) \left(x_{\gamma_k-\beta}^{\mathcal{M}_m}\right)^\intercal \right\rbrace \succeq 0
\end{align}

In this way, the relaxation order for the majority of buses in a large OPF problem can be set to a computationally tractable value with the computationally intensive higher-order relaxations only applied where necessary. 

\subsection{Iterative Solution Algorithm with a Heuristic for Determining the Relaxation Order}
\label{l:heuristic}

With a method for selectively applying the higher-order constraints to specific buses of an OPF problem, we next present an iterative solution algorithm for the moment relaxation. At each iteration, the algorithm uses a heuristic for specifying the value of $\gamma_k$ for each bus~$k$. Denote $\gamma$ as the vector containing $\gamma_k$, $k = 1,\dots,n$.

\begin{algorithm}
\caption{Iterative Solution for Moment Relaxation}\label{a:heuristic}
\begin{algorithmic}[1]
\State Set $\gamma_k = 1\quad \forall k \in \mathcal{N}$
\Repeat
\State Solve moment relaxation with $\gamma$
\State Calculate power injection mismatches
\State Increase the entries $\gamma$ according to the heuristic
\Until{Tolerances are satisfied}
\State Calculate optimal voltage profile
\end{algorithmic}
\end{algorithm}

Algorithm~\ref{a:heuristic} iteratively solves the moment relaxation and determines the power injection mismatches. The initial relaxation order $\gamma$ is set to one at every bus. If the solution meets specified tolerance criteria, the algorithm recovers the optimal voltage vector using an eigen decomposition of the diagonal block of the decomposed moment matrices corresponding to the second-order monomials. Otherwise, the relaxation order is increased at a subset of the buses with greatest power injection mismatch. Each iteration of the algorithm tightens the relaxation by adding higher-order constraints. The buses with greatest mismatch typically change with the relaxation order, thus potentially requiring multiple iterations of the loop in the algorithm.

There are several tolerance criteria used to evaluate the optimality of a solution. Due to numerical inaccuracies, no solver provides a solution that exactly satisfies the rank condition. One measure of the optimality of a candidate solution is based on power injection mismatches. Let $P^{mis}_k$ and $Q^{mis}_k$ be the active and reactive power mismatches at bus~$k$, respectively, resulting from the voltage vector derived from the closest matrix satisfying the rank condition. A voltage vector is accepted upon satisfaction of several convergence criteria: 1.) all apparent-power injection mismatches \mbox{$S_k^{mis} = \sqrt{\left(P^{mis}_k\right)^2 + \left(Q^{mis}_k\right)^2}$} are less than a specified tolerance, 2.) the voltage magnitudes, power injections, and line flows satisfy the inequality constraints in~\eqref{opf} to within specified tolerances, and 3.) the objective function evaluated with the voltage vector is equal to the optimal objective function from the moment relaxation to within a specified tolerance.

We finally describe the heuristic used to update $\gamma$. Define $\gamma^{\max}$ as the highest relaxation order among all buses (i.e., $\gamma^{\max} = \max_{k} \gamma_k$). (Note that $\gamma^{\max}$ is not a specified limit but rather can change as the algorithm progresses.) At each iteration of the algorithm, increment $\gamma_k$ at up to $h$ buses, where $h$ is a specified parameter, that have the largest apparent-power injection mismatches $S^{mis}_k$ among buses satisfying two conditions: 1.)~$\gamma_k$ is strictly less than $\gamma^{\max}$, and 2.)~$S_k^{mis}$ is greater than the specified tolerance. If no buses satisfy these two conditions, increment $\gamma_k$ at up to $h$ buses with the largest $S_k^{mis}$ greater than the specified tolerance and increment $\gamma^{\max}$. That is, in order to avoid unnecessarily increasing the size of the moment matrices, the heuristic avoids incrementing the maximum relaxation order $\gamma^{\max}$ until $\gamma_k = \gamma^{\max}$ at all buses $k$ with mismatch $S_k^{mis}$ greater than the mismatch tolerance.

There is a computational trade-off in choosing the value of $h$. Larger values of $h$ likely result in fewer iterations of the algorithm but each iteration is slower if more buses than necessary have high-order relaxations. Smaller values of $h$ result in faster solution at each iteration, but may require more iterations. Experience indicates that $h = 2$ is a good balance.

Note that this heuristic is just one of many possible approaches for specifying the relaxation order $\gamma_k$ at each bus~$k$. In addition to further analysis for varying $h$, future work includes comparison of this heuristic to alternative approaches.

\section{Numerical Results}
\label{l:results}

This section illustrates the effectiveness of the proposed algorithm for the higher-order relaxations relative to the first-order relaxation by considering several test problems. For these problems, second- and third-order relaxations are exact (i.e., the relaxation provides a solution with the globally minimal objective value and the globally optimal decision variables).

By modifying examples in the existing literature, it is straightforward to find moderate-size test problems for which existing solvers, such as the default interior point solver in \matpower~\cite{matpower}, fail for a variety of reasonable initialization heuristics. For instance, in an experiment conducted by randomly perturbing the cost function in the modified 118-bus system from~\cite{bukhsh_tps}, the default interior point solver in \matpower~with default solver options either fails to converge or converges to a local optimum in 6.9\% of 10000 tested problems initialized using five typical heuristics: 1.)~a ``flat start'' with unity voltage magnitudes and zero voltage angles, 2.)~the solution to the linear ``DC'' OPF approximation, 3.) a power flow solution calculated using active power injections at the midpoints of the generators' operating ranges, 4.)~a power flow solution calculated using power injections corresponding to an economic dispatch, and 5.)~a power flow solution calculated using the power injections resulting from a DC OPF. The first-order moment relaxation succeeds in globally solving all of the modified 118-bus test problems and other test problems for which traditional solution methods fail. These test problems are not of direct interest for our purposes but provide the context for considering the higher-order moment relaxations.

\begin{table}[tb]
\caption{Test Case Descriptions}
\label{t:testproblems}
\footnotesize
\renewcommand{\arraystretch}{1.3}
\centering
\begin{tabular}{|l|M{0.5\columnwidth}|c|}
\hline 
\multicolumn{1}{|c|}{\textbf{Test Case}} & \multicolumn{1}{c|}{\textbf{Description}} & \multicolumn{1}{|c|}{\textbf{Opt. Obj. Val.}}\\ 
& & \multicolumn{1}{|c|}{\textbf{(\$/hr)}} \\ \hline\hline
\emph{case14Q} & IEEE 14-bus system with all active and reactive loads decreased by 50\%\vspace{1pt} & $3.302\times 10^3$  \tabularnewline\hline
\emph{case14L} & IEEE 14-bus system with 25 MVA apparent-power flow limits for all lines\vspace{1pt} & $9.359\times 10^3$ \tabularnewline \hline
\emph{case39Q} & 39-bus system from~\cite{bukhsh_tps}, which is the IEEE 39-bus system with active and reactive loads decreased by 50\% and voltage bounds tightened to $\left[1.05, 0.95\right]$~per unit\vspace{1pt} & $1.122\times 10^4$ \tabularnewline \hline
\emph{case39L} & IEEE 39-bus system with all apparent-power line-flow limits decreased by 15\%\vspace{1pt} & $4.192\times 10^4$ \tabularnewline \hline
\emph{case57Q} & IEEE 57-bus system with active and reactive demand reduced by 75\% and all lower limits on generator reactive power injections set to -10~MVAr\vspace{1pt} & $7.352\times 10^3$ \tabularnewline \hline
\emph{case57L} & IEEE 57-bus system with 77 MVA apparent-power flow limits for all lines\vspace{1pt} & $4.398\times 10^4$ \tabularnewline \hline
\emph{case118Q} & IEEE 118-bus system with active and reactive demand reduced by 30\% and all lower limits on generator reactive power injection set to -20 MVAr\vspace{1pt} & $8.151\times 10^4$ \tabularnewline \hline
\emph{case118L} & IEEE 118-bus system with 110 MVA apparent-power flow limits for all lines\vspace{1pt} & $1.349\times 10^5$\tabularnewline \hline
\emph{case300} & IEEE 300-bus system\vspace{2pt} & $7.200\times 10^5$ \tabularnewline \hline
\multicolumn{3}{|M{0.95\columnwidth}|}{Small minimum resistances of $1\times 10^{-4}$~per unit are enforced on all branches in all test cases. All IEEE test cases are available in~\cite{ieee_test_cases}.\vspace{3pt}} \\\hline
\end{tabular}
\vspace{2pt}
\end{table}

\begin{table*}[t]
\caption{Test Case Results}
\label{t:results}
\small
\renewcommand{\arraystretch}{1.3}
\centering
\begin{tabular}{|c|l|c|c|c|c|c|c|}
\hline 
\textbf{Algorithm} & \multicolumn{1}{|c|}{\textbf{Test}} & \textbf{Max} $\mathbf{S^{mis}}$ & \textbf{Obj. Val.} & \textbf{Min. Eigenvalue} & \textbf{Num.} & \textbf{Solver Time} & \textbf{Num. High-Order} \\
 & \multicolumn{1}{|c|}{\textbf{Case}} & \textbf{(MVA)} & \textbf{Diff.} & \textbf{Ratio} & \textbf{Iter.} & \textbf{(sec)} & \textbf{Buses} \\ \hline\hline
\multirow{9}{0.18\columnwidth}{\centering Algorithm~\ref{a:heuristic}: Iterative Solver} 
 & \emph{case14Q}  & $1.08\times 10^{-3}$ & $2.36\times 10^{-6}$ & $1.08\times 10^6$ & 3 & 36.4 & ($2^{\mathrm{nd}}$): 3, ($3^{\mathrm{rd}}$): 0 \\ \cline{2-8}
 & \emph{case14L}  & $5.67\times 10^{-2}$ & $2.84\times 10^{-6}$ & $2.31\times 10^4$ & 3 & 25.5 & ($2^{\mathrm{nd}}$): 4, ($3^{\mathrm{rd}}$): 0  \\\cline{2-8} 
 & \emph{case39Q}  & $1.36\times 10^{-1}$ & $1.01\times 10^{-4}$ & $8.69\times 10^3$ & 19 & 2857 & ($2^{\mathrm{nd}}$): 31, ($3^{\mathrm{rd}}$): 2  \\\cline{2-8} 
 & \emph{case39L}  & $4.60\times 10^{-3}$ & $8.12\times 10^{-7}$ & $4.08\times 10^5$ & 2 & 4.88 & ($2^{\mathrm{nd}}$): 2, ($3^{\mathrm{rd}}$): 0  \\\cline{2-8}
 & \emph{case57Q}  & $6.49\times 10^{-3}$ & $4.35\times 10^{-6}$ & $2.52\times 10^5$ & 3 & 20.9 & ($2^{\mathrm{nd}}$): 4, ($3^{\mathrm{rd}}$): 0  \\ \cline{2-8}
 & \emph{case57L}  & $8.76\times 10^{-4}$ & $2.35\times 10^{-7}$ & $2.18\times 10^6$ & 2 & 88.8 & ($2^{\mathrm{nd}}$): 2, ($3^{\mathrm{rd}}$): 0  \\ \cline{2-8}
 & \emph{case118Q} & $2.13\times 10^{-1}$ & $3.40\times 10^{-5}$ & $6.78\times 10^4$ & 3 & 172.6 & ($2^{\mathrm{nd}}$): 4, ($3^{\mathrm{rd}}$): 0  \\ \cline{2-8}
 & \emph{case118L} & $4.42\times 10^{-1}$ & $4.59\times 10^{-5}$ & $1.93\times 10^4$ & 2 & 15.9 & ($2^{\mathrm{nd}}$): 2, ($3^{\mathrm{rd}}$): 0  \\ \cline{2-8}
 & \emph{case300}  & $5.14\times 10^{-2}$ & $3.74\times 10^{-6}$ & $4.65\times 10^4$ & 2 & 41.9 & ($2^{\mathrm{nd}}$): 2, ($3^{\mathrm{rd}}$): 0  \\ \hline\hline
\multirow{9}{0.18\columnwidth}{\centering Locally Minimal Set of Higher-Order Buses} 
 & \emph{case14Q}  & $7.64\times 10^{-3}$ & $3.21\times 10^{-5}$ & $1.52\times 10^5$ & N/A & 10.1 & ($2^{\mathrm{nd}}$): 2, ($3^{\mathrm{rd}}$): 0  \\ \cline{2-8}
 & \emph{case14L}  & $4.56\times 10^{-1}$ & $4.80\times 10^{-4}$ & $1.86\times 10^3$ & N/A & 8.51 & ($2^{\mathrm{nd}}$): 2, ($3^{\mathrm{rd}}$): 0  \\ \cline{2-8}
 & \emph{case39Q}  & $3.95\times 10^{-1}$ & $7.42\times 10^{-4}$ & $4.30\times 10^3$ & N/A & 289.7 & ($2^{\mathrm{nd}}$): 22, ($3^{\mathrm{rd}}$): 1  \\ \cline{2-8} 
 & \emph{case39L}  & $1.13\times 10^{-2}$ & $1.64\times 10^{-6}$ & $1.85\times 10^5$ & N/A & 3.21 & ($2^{\mathrm{nd}}$): 1, ($3^{\mathrm{rd}}$): 0  \\ \cline{2-8}
 & \emph{case57Q}  & $1.02\times 10^{-2}$ & $1.66\times 10^{-5}$ & $1.14\times 10^5$ & N/A & 7.67 & ($2^{\mathrm{nd}}$): 3, ($3^{\mathrm{rd}}$): 0  \\ \cline{2-8}
 & \emph{case57L}  & $1.34\times 10^{-3}$ & $6.60\times 10^{-7}$ & $6.15\times 10^5$ & N/A & 10.2 & ($2^{\mathrm{nd}}$): 1, ($3^{\mathrm{rd}}$): 0  \\ \cline{2-8}
 & \emph{case118Q} & \multicolumn{6}{c|}{Same as Algorithm~\ref{a:heuristic}}  \\ \cline{2-8}
 & \emph{case118L} & \multicolumn{6}{c|}{Same as Algorithm~\ref{a:heuristic}}  \\ \cline{2-8}
 & \emph{case300}  & \multicolumn{6}{c|}{Same as Algorithm~\ref{a:heuristic}}  \\ \hline\hline
 \multirow{9}{0.18\columnwidth}{\centering First-Order Relaxation}
 & \emph{case14Q}  & $4.92\times 10^{0}$  & $4.96\times 10^{-5}$ & $2.70\times 10^2$ & N/A & 0.69 & N/A \\ \cline{2-8}
 & \emph{case14L}  & $9.77\times 10^{0}$  & $5.94\times 10^{-4}$ & $1.41\times 10^2$ & N/A & 0.87 & N/A \\ \cline{2-8}
 & \emph{case39Q}  & $1.34\times 10^{2}$  & $3.54\times 10^{-2}$ & $2.89\times 10^2$ & N/A & 1.48 & N/A \\ \cline{2-8}
 & \emph{case39L}  & $6.91\times 10^{0}$  & $5.55\times 10^{-6}$ & $5.10\times 10^3$ & N/A & 1.41 & N/A \\ \cline{2-8}
 & \emph{case57Q}  & $5.76\times 10^{0}$  & $8.60\times 10^{-5}$ & $1.78\times 10^2$ & N/A & 2.08 & N/A \\ \cline{2-8}
 & \emph{case57L}  & $8.65\times 10^{0}$  & $1.58\times 10^{-3}$ & $3.84\times 10^2$ & N/A & 3.53 & N/A \\ \cline{2-8}
 & \emph{case118Q} & $7.36\times 10^{1}$  & $2.16\times 10^{-4}$ & $1.07\times 10^2$ & N/A & 4.81 & N/A \\ \cline{2-8}
 & \emph{case118L} & $1.07\times 10^{2}$  & $7.53\times 10^{-3}$ & $9.39\times 10^1$ & N/A & 7.64 & N/A \\ \cline{2-8}
 & \emph{case300}  & $2.42\times 10^{1}$  & $7.64\times 10^{-5}$ & $1.29\times 10^2$ & N/A & 17.9 & N/A \\ \hline
\end{tabular}
\vspace{2pt}
\end{table*}

Accordingly, it is also straightforward to modify examples in the existing literature to obtain test problems for which the first-order relaxation fails, but second- or third-order relaxations succeed. However, traditional solution methods, such as the interior point solver in \matpower, succeed in finding what turns out to be the global solution for these test problems. With the focus of this paper on demonstrating the effectiveness of higher-order moment relaxations relative to the first-order relaxation (and, equivalently, relative to the existing semidefinite relaxation~\cite{lavaei_tps}), we investigate these test problems. Thus, \matpower~with the default interior point solver and default solver options finds the global optimum for all the test problems in this section, with optimal objective values listed in Table~\ref{t:testproblems} and solution times of less than 0.5 seconds, but the first-order relaxation fails to solve each test problem. Second- and third-order moment relaxations \emph{certify} that these solutions are the global optima.

Note that the existence of small examples for which both the first-order relaxation and traditional solution methods fail but the higher-order moment relaxations succeed in finding the global optimum (e.g., the five- and nine-bus systems in~\cite{bukhsh_tps} and the five-bus system in~\cite{bernie_opfconvexity}; see~\cite{pscc2014,patrick,ibm_paper} for analysis of these and other small test cases) implies that similar phenomena can occur in large practical problems as well. This suggests the need for a wider variety of test problems, the development of which is beyond the scope of this paper.

Descriptions of each test problem are provided in Table~\ref{t:testproblems}. Except for the IEEE 300-bus system, these problems are modifications of the IEEE test cases using one of two known methods for inducing failure of the first-order relaxation. Tightening apparent-power line-flow limits may induce failure of the first-order relaxation~\cite{hicss2014}. The letter ``L'' denotes corresponding problems. Decreasing the loading while reducing the generators' leading power factor range (i.e., decreasing the magnitude of the lower reactive power generation limits) may also induce failure of the first-order relaxation~\cite{bukhsh_tps}. The letter ``Q'' denotes corresponding problems.

The results in this section are generated using a computer with a quad-core 2.70 GHz processor and 16 GB of RAM. The moment relaxations are implemented using \mbox{MATLAB 2013a}, YALMIP version 2014.02.21~\cite{yalmip}, and Mosek version 7.0.0.102~\cite{mosek}.

Table~\ref{t:results} presents the results from applying the moment relaxations to the test problems in Table~\ref{t:testproblems}. The first group of rows in Table~\ref{t:results} shows the results from Algorithm~\ref{a:heuristic}. The second group shows the results from an at-least-locally minimal set of higher-order buses. This set is derived by individually removing higher-order buses from the solution given by Algorithm~\ref{a:heuristic}. The third group shows the results from the first-order relaxation.

The columns show the values of three convergence metrics. The first metric is the maximum apparent-power injection mismatch (Max $S^{mis}$), which has a 0.5 MVA tolerance. Note that the voltage magnitudes, power injections, and line flows satisfy the inequality constraints in~\eqref{opf} to within 0.005~per~unit voltage and 0.5 MVA for both Algorithm~\ref{a:heuristic} and the locally minimal set of higher-order buses.

The second metric compares the optimal objective value obtained directly from the moment relaxation, denoted as $C^{\mathrm{Mom.\,Obj.}}$ and the cost implied by the voltage vector obtained from the closest matrix satisfying the rank condition, denoted as $C^{V}$. In the absence of numerical inaccuracy, these objective values are equal when the solution to the moment relaxation is exact. With imperfect solvers, this value may be non-zero even when the moment relaxation has an exact solution. The column Obj. Val. Diff. shows $\left|C^{\mathrm{Mom.\,Obj.}} - C^{V} \right|/C^{\mathrm{Mom.\,Obj.}}$, which is employed as the second convergence criterion with a tolerance of $1\times 10^{-3}$. For a solution with small power injection mismatches, a small value in this column indicates global optimality for practical purposes. 

\begin{figure*}[!t]
\centering
\includegraphics[totalheight=0.52\textheight]{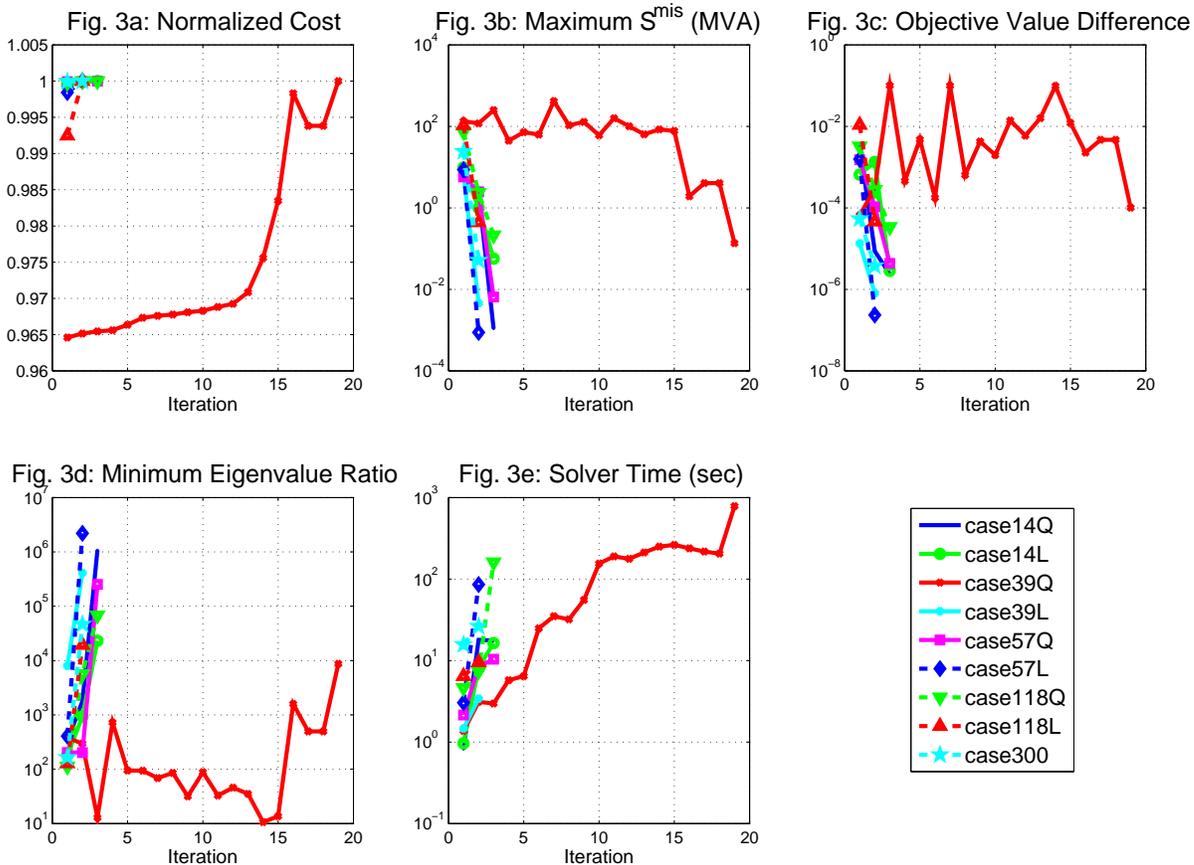}
\caption{Test Case Results for Each Iteration of Algorithm~\ref{a:heuristic}. Fig.~\ref{f:quality_and_time}a shows the lower bound on the optimal objective value normalized by the globally optimal objective value. Figs.~\ref{f:quality_and_time}b, ~\ref{f:quality_and_time}c, and ~\ref{f:quality_and_time}d show the three convergence metrics, and Fig.~\ref{f:quality_and_time}e shows the solver time for each iteration.}\label{f:quality_and_time} 
\end{figure*}

Since the first-order relaxations yield large power injection mismatches, the objective value calculated from the voltage vector $\left(C^{V} \right)$ is not meaningful. Therefore, the Obj. Val. Diff. column for the first-order relaxation shows the difference between the globally optimal objective value $C^{\mathrm{Mom.\,Obj.}}$ and the lower bound obtained from the objective value of the first-order relaxation $C^{\mathrm{SDP}}$ as $\left(C^{\mathrm{Mom.\,Obj.}} - C^{\mathrm{SDP}} \right)/C^{\mathrm{Mom.\,Obj.}}$.

The third metric, which is first proposed in~\cite{molzahn_holzer_lesieutre_demarco-large_scale_sdp_opf}, is based directly on satisfaction of the rank condition. For a solution that satisfies the rank condition, the diagonal blocks of the moment matrices corresponding to the second-order terms have rank equal to one (i.e., the matrices $\mathbf{M}_1^{\mathcal{M}_m}\left(y\right)$ for $m=1,\ldots,\left|\mathcal{M}\right|$ have a single non-zero eigenvalue). However, for numerical reasons, solvers do not yield a ``hard zero'' value for ``zero'' eigenvalues of these matrices. On the other hand, when the moment relaxation fails to be exact, these matrices have more than one non-zero eigenvalue.

To measure the satisfaction of the rank condition, we use the ratio between the largest and second-largest magnitude eigenvalues. The minimum such ratio among all the matrices $\mathbf{M}_1^{\mathcal{M}_m}\left(y\right)$ for $m=1,\ldots,\left|\mathcal{M}\right|$, is termed the \emph{minimum eigenvalue ratio}. If the solution to a case for which the relaxation was exact did have ``hard zeros'' for zero eigenvalues, the largest eigenvalue would be non-zero and the second-largest eigenvalue would be zero, resulting in a minimum eigenvalue ratio of infinity. In practice, numerical issues result in minimum eigenvalue ratios that are large (typical values are on the order of $10^4$ for problems that satisfy the rank condition). Further, if the solution does not satisfy the rank condition, both the largest and second largest eigenvalues typically have similar magnitudes, therefore yielding a small value for the minimum eigenvalue ratio. Thus, a large value for the minimum eigenvalue ratio indicates satisfaction of the rank condition while a small value indicates failure to satisfy the rank condition.

Note that we do not use the minimum eigenvalue ratio as a convergence criteria. This is due to the fact that a solution to the moment relaxation can have a relatively poor (small) minimum eigenvalue ratio but the closest rank one matrix can still globally solve the OPF problem to within the other tolerances. (See, for instance, the results for \emph{case14L} and \emph{case39Q} in the second group of rows in Table~\ref{t:results} which show the results for the locally minimal set of higher-order buses.) We report the minimum eigenvalue ratio in the column Min. Eigenvalue Ratio of Table~\ref{t:results}.

Where applicable, Table~\ref{t:results} also shows the number of iterations of Algorithm~\ref{a:heuristic} (Num. Iter.), the time spent in the MOSEK solver summed over all iterations (Solver Time), and the number of higher-order buses (Num. High-Order Bus), where the quantity in parentheses is the relaxation order. 

Fig.~\ref{f:quality_and_time} shows the lower bound on the objective function, the convergence metrics, and the solver times for each iteration of Algorithm~\ref{a:heuristic}. The lower bounds on the objective function in Fig.~\ref{f:quality_and_time}a are normalized so that the global optimum for each test problem has a value of one. Using a log scale, Fig.~\ref{f:quality_and_time}b shows the maximum power injection mismatch in MVA (i.e., the first convergence metric) for each iteration, Fig.~\ref{f:quality_and_time}c shows the objective value difference (i.e., the second convergence metric), and Fig.~\ref{f:quality_and_time}d shows the minimum eigenvalue ratio (i.e., the third convergence metric). Fig.~\ref{f:quality_and_time}e shows the solver time for each iteration of the algorithm.

We next emphasize several interesting observations from the results in Table~\ref{t:results} and Fig.~\ref{f:quality_and_time}. Solutions to the first-order relaxations generally have a subset of buses with large power injection mismatches, but the objective values are often quite close to the global optimum (see column Obj. Val. Diff. in Table~\ref{t:results} and the fact that the first iteration in Fig.~\ref{f:quality_and_time}a is approximately equal to 1 but the convergence metrics in Figs.~\ref{f:quality_and_time}b, \ref{f:quality_and_time}c, and \ref{f:quality_and_time}d are not satisfied). This suggests that there are often ``hidden'' or ``nearly hidden'' rank one solutions for the first-order relaxation~\cite{bitar_allerton2013}. Selective application of the higher-order constraints provides a mechanism for recovering these hidden rank one solutions. However, not all problems have hidden rank one solutions (e.g., the first-order relaxations of \emph{case39Q} and \emph{case118L}, which yield solutions that are 3.54\% and 0.75\%, respectively, below the global optimum, and the examples in~\cite{pscc2014}). The solver times from Algorithm~\ref{a:heuristic} illustrate that additional computational effort is required to achieve the global solution relative to the lower bounds from the first-order relaxation.

The Num. Iter. column in Table~\ref{t:results} and the convergence metrics in Figs.~\ref{f:quality_and_time}b, \ref{f:quality_and_time}c, and \ref{f:quality_and_time}d show that only a small number of iterations are typically required to obtain a global optimum (i.e., the higher-order moment constraints are only required at a small number of buses). However, this is not always the case as \emph{case39Q} requires many iterations and higher-order constraints at the majority of buses in the network. This problem demonstrates that the approach of selectively applying higher-order constraints to specific areas of the network may not be computationally tractable for all problems.

Since each iteration of Algorithm~\ref{a:heuristic} adds constraints to the optimization problem, the cost shown in Fig.~\ref{f:quality_and_time}a should be non-decreasing. The cost is non-decreasing for all problems with the exception of iteration 17 for \emph{case39Q}. The decrease at iteration 17 is explained by the fact that the semidefinite programming solver does not converge to a sufficient tolerance (i.e., the constraints are not satisfied) at this iteration of Algorithm~\ref{a:heuristic}.

For \emph{case118Q}, \emph{case118L}, and \emph{case300}, Algorithm~\ref{a:heuristic} finds an at-least-locally minimal set of buses requiring the higher-order relaxation. For \emph{case14Q}, \emph{case39L}, \emph{case57Q}, and \emph{case57L}, Algorithm~\ref{a:heuristic} only uses one more bus than a locally minimal set. While this indicates that the heuristic with $h=2$ is effective in identifying a minimal or near-minimal set of higher-order buses for many OPF problems, we emphasize that these are only known to be locally minimal; a different heuristic may identify a smaller set of buses. Further, Algorithm~\ref{a:heuristic} does not identify a near-minimal set of higher-order buses for \emph{case39Q}, which results in an almost order-of-magnitude larger solution time than necessary. More sophisticated heuristics could lead to better performance for some problems.



Finally, note that solution times have stronger dependence on the number of higher-order buses, the size of the maximal cliques they are contained within, and the relaxation orders required than the size of the system. For instance, although \emph{case300} has 7.7 times more buses, Algorithm~\ref{a:heuristic} has a factor of 68.2 greater solution time for \emph{case39Q} due to the number of buses with second- and third-order constraints in \emph{case39Q}. Further note that although solution times for the moment relaxations are not yet competitive with mature local solvers, such as the interior point method in \matpower~which solved all test problems in Table~\ref{t:testproblems} in less than 0.5 seconds, the moment relaxations with Algorithm~\ref{a:heuristic} provide a computationally tractable approach for globally solving many problems for which the first-order relaxation fails to yield a global solution. The moment relaxations also \emph{certify} global optimality in contrast to traditional solvers which only guarantee a local optimum.

\section{Conclusion}
\label{l:conclusion}

While existing convex relaxations globally solve many OPF problems, there are practical problems for which existing relaxations fail to yield physically meaningful solutions. This paper has described a hierarchy of ``moment''  relaxations that globally solve many problems for which existing relaxations fail. The moment relaxations, which take the form of semidefinite programs, are developed from the Lasserre hierarchy for generalized moment problems. Increasing the order in this hierarchy results in ``tighter'' relaxations at the computational cost of larger semidefinite programs.

Solving the moment relaxations for larger problems requires both exploiting power system sparsity and selectively applying the higher-order moment relaxation constraints. A matrix completion decomposition for exploiting sparsity was first presented. Next, taking advantage of the observation that first-order relaxations are sufficient for large regions of typical OPF problems, this paper proposed an iterative algorithm for solving the moment relaxations. A heuristic at each iteration identifies where to enforce the higher-order relaxation constraints. The proposed algorithm's effectiveness was demonstrated by globally solving several test cases for which existing convex relaxations failed.

Future work includes improving the computational performance of the proposed algorithm. Alternative heuristics for determining where to apply the higher-order constraints may reduce solution times. Distributed solution algorithms, which have proven valuable for existing relaxations~\cite{lam2012}, may also speed computation of the moment relaxations.

Identification and exploration of realistic test cases for which low-order relaxations fail is another important future research direction. The problems in~\cite{lavaei_tps} used to demonstrate that the OPF problem is, in general, NP-hard serve as example cases for which low-order moment relaxations fail to yield a global solution~\cite{pscc2014}. However, with a large number of global optima, these problems are very atypical. Extending the work of, e.g.,~\cite{low_tutorial1,low_tutorial2,madani2014}, sufficient conditions for tightness of the moment relaxations would also be valuable contributions. 

Additional future work also includes application of moment relaxations to more general OPF formulations that include, for instance, discrete devices, security constraints, and transient stability constraints. Exploiting synergies with robust and chance-constrained optimization techniques for the OPF problem appears particularly promising~\cite{pscc2014survey}. Extension to other problems in power system engineering, such as state estimation, voltage stability margins, and power flow calculations, is another avenue of future work.

\section*{Acknowledgment}
The authors acknowledge the support of the Dow Postdoctoral Fellowship in Sustainability, ARPA-E grant \mbox{DE-AR0000232}, and Los Alamos National Laboratory subcontract 270958.

\vfill
\pagebreak
\begin{IEEEbiography}[{\includegraphics[width=1in,height=1.25in,clip,keepaspectratio]{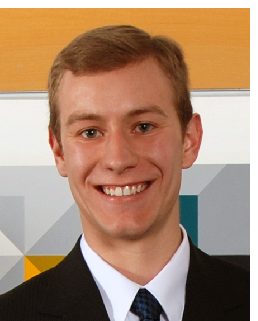}}]{Daniel K. Molzahn}
(S'09-M'13) is a Dow Sustainability Fellow at the University of Michigan, Ann Arbor. He received the B.S., M.S., and Ph.D. degrees in electrical engineering and the Masters of Public Affairs degree from the University of Wisconsin-–Madison, where he was a National Science Foundation Graduate Research Fellow. His research interests are in the application of optimization techniques and policy analysis to electric power systems.
\end{IEEEbiography}

\begin{IEEEbiography}[{\includegraphics[width=1in,height=1.25in,clip,keepaspectratio]{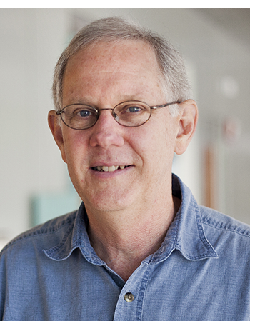}}]{Ian A. Hiskens} (S'77-M'80-SM'96-F'06) received the B.Eng. degree in electrical engineering and the B.App.Sc. degree in mathematics from the Capricornia Institute of Advanced Education, Rockhampton, Australia, in 1980 and 1983 respectively, and the Ph.D. degree in electrical engineering from the University of Newcastle, Australia, in 1991.

He is the Vennema Professor of Engineering in the Department of Electrical Engineering and Computer Science, University of Michigan, Ann Arbor. He has held prior appointments in the Queensland electricity supply industry, and various universities in Australia and the United States. His research interests lie at the intersection of power system analysis and systems theory, with recent activity focused largely on integration of renewable generation and controllable loads.

Dr. Hiskens is actively involved in various IEEE societies, and is VP-Finance of the IEEE Systems Council. He is a Fellow of the IEEE, a Fellow of Engineers Australia, and a Chartered Professional Engineer in Australia.
\end{IEEEbiography}

\vfill


\end{document}